\documentclass[12pt,leqno]{article}
\usepackage{amsmath, amsthm, amsfonts, amssymb, color}
\usepackage{mathrsfs,enumerate}
\newtheorem{thm}{Theorem}[section]
\newtheorem{cor}[thm]{Corollary}

\newtheorem{prp}[thm]{Proposition}

\theoremstyle{definition}
\newtheorem{defn}{Definition}[section]
\newcommand{\scr}[1]{\mathscr #1}
\definecolor{wco}{rgb}{0.5,0.2,0.3}

\numberwithin{equation}{section} \theoremstyle{remark}
\newtheorem{rem}{Remark}[section]

\newcommand{\ua}{\uparrow}

\title{{\bf Strong Solutions of  Stochastic Generalized Porous Media Equations:
Existence, Uniqueness and Ergodicity}
\footnote{Supported in part by the DFG through the Forschergruppe
``Spectral Analysis, Asymptotic
Distributions and Stochastic Dynamics'', the BiBoS Research Centre,
NNSFC(10025105,10121101), TRAPOYT
in China and the 973-Project.} }
\author{{\bf G. Da Prato$^{a)}$, Michael R\"ockner$^{b)}$, B.L. Rozovskii$^{c)}$ and Feng-Yu Wang$^{d)}$}\\
\footnotesize{a) Scuola Normale Superiore, Piazza dei Cavalieri 7, 56126, Pisa, Italy}\\
\footnotesize{b) Fakult\"at f\"ur Mathematik, Universit\"at Bielefeld, D-33501 Bielefeld, Germany}\\
\footnotesize{c) Center of Applied Mathematics, USC, Los Angeles, CA, 90089-1113 USA}\\
\footnotesize{d) School of Mathematical Science, Beijing Normal University, Beijing 100875, China}}
\begin{document}
\maketitle
\begin{abstract} Explicit conditions are presented for the existence, uniqueness and
ergodicity of the strong solution to a class of generalized stochastic porous media
equations. Our estimate of the convergence rate is sharp
according to the known optimal decay for the solution of the classical
(deterministic) porous medium equation.

\end{abstract}
 \vskip 2cm

\def\C{\scr C}
\def\D{\scr D}
\def\E{\scr E}
\def\F{\scr F}

\def\BB{\mathbb B}
\def\EE{\mathbb E}
\def\N{\mathbb N}
\def\R{\mathbb R}
\def\Z{\mathbb Z}

\def\d{\text{\rm{d}}}
\def\m{{\bf m}}

\def\aa{\alpha}
\def\bb{\beta}
\def\dd{\delta} \def\DD{\Delta}
  \def\ee{\varepsilon}
  \def\vv{\varepsilon}
\def\gg{\gamma} \def\GG{\Gamma}
\def\kk{\kappa}
\def\ll{\lambda}
\def\rr{\rho}
  \def\vrr{\varrho}
\def\si{\sigma}
\def\oo{\omega}  \def\OO{\Omega}

\def\ff{\frac} \def\ss{\sqrt}
\def\<{\langle} \def\>{\rangle}
  \def\nn{\nabla} \def\pp{\partial}

\def\ess{\text{\rm{ess}}}
\def\beg{\begin} \def\beq{\begin{equation}}
\def\Ric{\text{\rm{Ric}}} \def\Hess{\text{\rm{Hess}}}\def\B{\mathbb B}
\def\e{\text{\rm{e}}} \def\ua{\underline a}  \def\b{\mathbf b}
\def\tt{\tilde}
\def\cut{\text{\rm{cut}}} \def\P{\mathbb P} \def\ifn{I_n(f^{\bigotimes n})}
\def\fff{f(x_1)\dots f(x_n)} \def\ifm{I_m(g^{\bigotimes m})}
\def\pm{\pi_{{\bf m}}}   \def\p{\mathbf{p}}   \def\ml{\mathbf{L}}
    \def\aaa{\mathbf{r}}     \def\r{r}
\def\gap{\text{\rm{gap}}} \def\prr{\pi_{{\bf m},\varrho}}  \def\r{\mathbf r}

\section{Main Results}

Let $(E,\scr M,{\bf m})$ be a separable probability space and  $(L,\D(L))$ a
negative definite self-adjoint linear operator on $L^2({\bf m})$ having
 discrete spectrum with eigenvalues
\[
   0 > -\ll_1\ge -\ll_2\ge \cdots\to -\infty
\]
and $L^2(m)$-normalized eigenfunctions $\{e_i\}$ such that $e_i \in L^{r+1}({\bf
m})$ for any $i\ge 1$, where $r>1$ is a fixed number
throughout this paper. We assume
that $L^{-1}$ is bounded in $L^{r+1}(\m)$, which is e.g.\ the case if
$L$ is a Dirichlet operator
(cf.\ e.g.\ \cite{ma-roe92})
since in this case the interpolation
theorem
or simply Jensen's inequality
implies $\|\e^{tL}\|_{r+1}\le \e^{-\ll_1t 2/(r+1)}$ for
all $t\ge 0.$ A classical example of $L$ is the Laplacian operator
on a smooth bounded domain in a complete Riemannian manifold with
Dirichlet boundary conditions.

In this paper we consider the following stochastic differential equation:
\beq\label{1.2}
  \d X_t = (L\Psi(X_t)+\Phi(X_t))\d t+Q\d W_t,
\end{equation}
where $\Psi$ and $\Phi$ are (non-linear)  continuous functions on $\R$,
and $Q$ is a densely defined linear operator on $L^2({\bf m})$
with $Q e_i:= \sum_{j=1}^\infty q_{ji}  e_j\ (i\ge 1)$ such that
$q_i^2:=\sum_{j=1}^\infty q_{ij}^2$ satisfies
\[
  q:=\sum_{i=1}^\infty \ff{q_i^2}{\ll_i}<\infty.
\]
An appropriate Hilbert space $H$ as state space for the solutions
to (\ref{1.2}) is given as follows. Let
\[
  H^1 := \Bigl\{
           f\in L^2(m) \,:\, \sum_{i=1}^\infty \lambda_i \m(f e_i)^2 <\infty
         \Bigr\}.
\]
Define $H$ to be its topological dual with inner product $\< \, , \, \>_H$.
Identifying  $L^2(\m)$ with its dual we get the continuous and dense
embeddings
\[
  H^1 \subset L^2(\m) \subset H.
\]
We denote the duality between $H$ and $H^1$ by $\<\,,\,\>$.
Obviously, when restricted to $L^2(\m)\times H^1$ this coincides
with the natural inner product in $L^2(\m)$, which we therefore also
denote by $\<\,,\,\>$, and it is also clear that
\[
  \<f,g\>_H
  = \sum \lambda_i^{-1} \<f,e_i\> \<g,e_i\>,
  \quad f,g\in H.
\]
Furthermore, in (\ref{1.2})
$W_t = (b_t^i)_{i\in\N}$ is a cylindrical Brownian motion on $L^2(\m)$
where $\{b_t^i\}$ are independent one-dimensional Brownian motions
on a complete probability space $(\OO, \F, P)$. Let  $\F_t$ be the natural
filtration of $W_t$. Then
\[
  Q W_t
  := \sum_{i=1}^\infty \Bigl(\sum_{j=1}^\infty q_{ij} b_t^{j}\Bigr)e_i,
  \quad t\geq 0 ,
\]
is a well-defined process taking values in $H$ which is a martingale.

Recall that the classical porous medium equation reads
\[
  \d X_t = \DD X_t^m \d t
\]
on a domain in $\R^d$, see e.g. \cite{Aronson} and the
references therein. 
So, we may call (\ref{1.2}) the generalized stochastic porous medium
equation.
Recently, the existence and uniqueness of
weak solutions as well as the existence of invariant probability measures
for the stochastic porous medium equation, i.e. (\ref{1.2}) with $L:=\DD$
on a bounded domain in $\R^d$ with Dirichlet
boundary conditions, were proved in \cite{DR1, DR2, BDR, BBDR}.
In this paper we aim to prove the existence and uniqueness of strong solutions
of (\ref{1.2}), in particular, describe the convergence rate of the solution,
for $t\to\infty$.

To solve (\ref{1.2}), we assume that there exist some constants
$c \geq 0$, $\eta,\si\in\R$ such that
\beq\label{1.3}
  \beg{split}
    & |\Psi(s)| +|\Phi(s)| \le c (1+|s|^r),\\
    & (s-t)(\Psi(s)-\Psi(t))\ge \eta |s-t|^{r+1}+\si (s-t)^2,
      \ \ \ s,t\in \R.
  \end{split}
\end{equation}
Since by the Cauchy-Schwarz inequality one has
\beg{equation*}
  \beg{split}
    \ff{2^{3-r}|s-t|^{r+1}}{(r+1)^2}
      &\le\ff{4}{(r+1)^2}\big(|s|^{(r+1)/2}\mathop{\text{sgn}}(s)
       - |t|^{(r+1)/2}\mathop{\text{sgn}}(t)\big)^2\\
      &=\bigg(\int_t^s |u|^{(r-1)/2}\d u\bigg)^2 \le (s-t)\int_t^s |u|^{r-1}\d u,
  \end{split}
\end{equation*}
(\ref{1.3}) holds if $\Psi(0)=0$ and there exists $\kk>0$ such that (cf. \cite{BDR})
\[
  \si +\ff{(r+1)^2} 4 \eta |s|^{r-1}\le \Psi'(s) \le \kk (1+|s|^{r-1}),\ \ \ s\in\R.
\]
Next, assume that there exist $\theta< \eta$ and $\dd\le \si$
such that
\beq\label{1.3'}
  -{\bf m} \big((\Phi(x)-\Phi(y))L^{-1} (x-y)\big)
  \le \theta \|x-y\|_{r+1}^{r+1} +\dd \|x-y\|_2^2,
  \ \ \ x,y\in L^{r+1}({\bf m}),
\end{equation}
where here and in the sequel, $\|\cdot\|_p$ denotes the $L^p$-norm
with respect to ${\bf m}$ for any $p\ge 1.$ We note that since
$L^{-1}$ is bounded on $L^{r+1}(\m)$ and $r>1$, if there exist
constants $c_1,c_2\ge 0$ such that

$$|\Phi(s)-\Phi(t)|\le c_1 |s-t|^r +c_2 |s-t|,\ \ \ s,t\in\R,$$
then

$$-{\bf m}\big((\Phi(x)-\Phi(y)) L^{-1}(x-y)\big) \le c_1
\|L^{-1}\|_{r+1}\|x-y\|_{r+1}^{r+1}+ c_2\ll_1^{-1} \|x-y\|_2^2,$$
hence (\ref{1.3'}) holds for $\theta:= c_1\|L^{-1}\|_{r+1}$ and
$\dd:= c_2\ll_1^{-1}.$

%

\begin{defn}\label{defn1.1}
  Let $\nu(\d t):= \e^{-t} \d t.$
  An $H$-valued continuous $(\scr F_t)$-adapted process $X_t$
  is called a solution to (\ref{1.2}),
  if $X\in L^{r+1}(\R_+\times \OO\times E; \nu\times P\times \m)$
  such that
  \beq\label{D}
    \<X_t,e_i\>
    = \<X_0, e_i\>
      + \int_0^t {\bf m}\big(\Psi(X_s)Le_i+ \Phi(X_s) e_i\big) \d s +q_i B_t^i,
      \ \ i\ge 1, t>0,
  \end{equation}
  where $B_t^i := \ff 1 {q_i} \sum_{j=1}^\infty q_{ij} b_t^j$
  ($:= 0$ if $q_i=0$) is an $(\F_t)$-Brownian motion on $\R$
  (provided it is non-trivial).
\end{defn}

%
%

\begin{rem}\label{rem1.1}
  \begin{enumerate}[(i)]
    \item We note that (\ref{D}) indeed makes sense, since by the first
      inequality in (\ref{1.3}) we have
      \[
        \Psi(X), \Phi(X)
        \in L^{(r+1)/r} (\R_+ \times \Omega \times E, \, \nu\times P\times \m).
      \]
    \item We emphasize that for each solution of (\ref{D}) there also
      exists a vector-valued version of the equation. More precisely, the
      integral comes from an $H$-valued random vector with a natural integrand
      which, however, takes values in a larger Banach space $\BB'$.
      To describe this in detail, we need some preparations.
  \end{enumerate}
\end{rem}
Consider the separable Banach space $\BB := L^{r+1}(\m)$.
Then we can obtain a presentation of its dual space $\BB'$
through the embeddings
\[
  \BB \subset H \equiv H' \subset \BB' ,
\]
where $H$ is identified with its dual through the Riesz-isomorphism.
In other words $\BB'$ is just the completion of $H$ with respect to
the norm
\[
  \|f\|_{\BB'}
  := \sup_{\|g\|_{r+1}\leq 1} \<f,g\>_H,
  \quad f\in H.
\]
Since $H$ is separable, so is $\BB'$.
We note that this is different from the usual representation of $\BB=L^{r+1}(\m)$
through the embedding
\[
  \BB \subset L^2(\m) \equiv L^2(\m)' ,
\]
which, of course, gives $L^{(r+1) / r}(\m)$ as dual. But it is
easy to identify the isomorphism between $L^{(r+1) / r}(\m)$ and
$\BB'$. Below $_{\BB'}\<\,,\,\>_{\BB}$ denotes the duality
between $\BB$ and $\BB'$. Clearly, $_{\BB'}\<\,,\,\>_\BB =
\<\,,\,\>_H$ on $\BB\times H$.
\begin{prp}\label{prp1.1}
  The linear operator
  \[
    Lf
    := -\sum_{i=1}^\infty \lambda_i \, \m(f e_i)e_i,
    \quad f\in L^2(\m),
  \]
  defines an isometry from $L^{(r+1)/r}(\m)$ to $\BB'$
  with dense domain.
  Its (unique) continuous extension $\bar L$
  to all of $L^{(r+1) / r}(\m)$ is an isometric
  isomorphism from $L^{(r+1)/r}(\m)$ onto $\BB'$ such that
  \begin{equation}\label{eq1.5}
    _{\BB'}\< -\bar L f, g \>_{\BB}
    = \m(fg)
    \quad\text{for all }f\in L^{(r+1)/r}(\m), \; g\in L^{r+1}(\m).
  \end{equation}
\end{prp}
\begin{proof}
  Let $f\in L^2(\m)$, $N>n\geq 1$. Then
  \begin{align*}
    \Bigl\| \sum_{i=n}^N \lambda_i\, \m(f e_i)e_i \Bigr\|_{\BB'}
    &= \sup_{\|g\|_{r+1}\leq 1} \Bigl| \m\Bigl( g\sum_{i=n}^N \m(f e_i)e_i \Bigr) \Bigr| \\
    &= \Bigl\| \sum_{i=n}^N \m(f e_i)e_i \Bigr\|_{(r+1)/r} .
  \end{align*}
  Since $f=\sum_{i=1}^\infty \m(f e_i)e_i$ with the series converging
  in $L^2(\m)$, hence in $L^{\frac{r+1}{r}}(\m)$ (because $r>1$),
  the first part of the assertion follows, and $\bar L$ is an isometry
  from $L^{(r+1)/r}(\m)$ into $\BB'$.
  Now let $T\in\BB'$. Then there exists $f\in L^{(r+1)/r}(\m)$
  such that for all $g\in L^{r+1}(\m)$
  \begin{equation}\label{eq1.6}
    \begin{split}
      _{\BB'}\<T,g\>_{\BB}
      &= \m(fg)\\
      &= \lim\limits_{n\to\infty}\m(f_n g)
    \end{split}
  \end{equation}
  for some $f_n\in D(L)$ such that $\lim\limits_{n\to\infty}f_n =f$
  in $L^{(r+1)/r}(\m)$. Hence for all $g\in L^{r+1}(\m)$
  \begin{equation}\label{eq1.7}
    \begin{split}
      _{\BB'}\<T,g\>_{\BB}
      &= \lim_{n\to\infty} \m\bigl(f_n \, L(L^{-1}g)\bigr)\\
      &= \lim_{n\to\infty} \m(Lf_n \, L^{-1}g) \\
      &= \lim_{n\to\infty} \<-Lf_n, g\>_H
        = - \, _{\BB'}\< \bar Lf , g \>_{\BB}
    \end{split}
  \end{equation}
  and the second assertion is proved.
  Since any $f\in L^{(r+1)/r}(\m)$ defines a $T\in\BB'$, the last
  assertion follows from (\ref{eq1.6}) and (\ref{eq1.7}).
\end{proof}
Since $L^{-1}$ is bounded on $L^{r+1}(\m)$ by our assumptions on $L$
(which was not used so far), we obtain the following consequence.
\begin{cor}\label{cor1.1}
  Let $(L^{-1})' : L^{(r+1)/r}(\m) \to L^{(r+1)/r}(\m)$ be the dual
  operator of $L^{-1}: L^{r+1}(\m)\to L^{r+1}(\m)$.
  Then the operator
  \[
    J : \bar L \circ (L^{-1})' : L^{(r+1)/r}(\m) \to \BB'
  \]
  extends the natural inclusion $L^2(\m)\subset H\subset \BB'$
  and for all $f\in L^{(r+1)/r}(\m)$
  \beq\label{eq1.7.1}
    _{\BB'}\< Jf , g \>_{\BB}
    = -\m(f\, L^{-1}g)
    \quad\text{for all }g\in L^{r+1}(\m) .
  \end{equation}
\end{cor}
\begin{proof}
  If $f\in L^2(\m)$, then $(L^{-1})'\, f = L^{-1}f$,
  hence $Jf = f \in L^2(\m)\subset H\subset \BB'$.
  The last assertion follows by (\ref{eq1.5}).
\end{proof}
Multiplied by $\lambda_i^{-1}$, (\ref{D}) by the above now reads
\[
  _{\BB'}\< X_t , e_i \>_{\BB}
  = {} _{\BB'}\< X_0 , e_i \>_{\BB}
    + \int_0^t {}_{\BB'}\bigl\< \bar L \Psi(X_s) + J\Phi(X_s) , \, e_i \bigr\>_{\BB} \,\d s
    + {} _{\BB'}\< QW , e_i \>_{\BB}.
\]
By Remark~\ref{rem1.1}, Proposition~\ref{prp1.1} and Corollary~\ref{cor1.1}
the Bochner integrals
\[
  \int_0^t \bigl(\bar L \Psi(X_s) + J \Phi(X_s)\bigr) \,\d s,
  \quad t\geq 0,
\]
exist in $\BB'$. So, (\ref{D}) can always be rewritten
equivalently in vector form as an
equation in $\BB'$ as
\begin{equation}\label{eq1.8}
  X_t
  = X_0 + \int_0^t \bigl( \bar L \Psi(X_s) + J\Phi(X_s) \bigr)\,\d s
    + QW_t,
  \quad t\geq 0.
\end{equation}
Note that by Definition~\ref{defn1.1}, $X_t\in H$ and also $QW_t \in H$,
hence the integral in (\ref{eq1.8}) is necessarily a continuous
$H$-valued process.

Now we can state our main results.
\beg{thm}\label{T1.1}
  Assume $(\ref{1.3})$ and $(\ref{1.3'})$  with $\si\ge \dd$ and $ \eta> \theta$.
  We have:
  \begin{enumerate}[(1)]
    \item For any $\F_0 / \scr B(H)$-measurable $\xi : \Omega\to H$
      with $\EE \|\xi\|_H^2 <\infty$ there exists a unique solution
      $X$ to (\ref{1.2}) such that $X_0 = \xi$.
      Furthermore, there exists $C>0$ such that
      \beq\label{eq1.8.1}
        \EE \| X_t \|_H^2 \leq C (1+t^{-2/(r-1)}),
        \quad t>0.
      \end{equation}
      In particular,
      for any $x\in H$ there exists a unique solution  $X_t(x)$ to $(\ref{1.2})$
      with initial value $x$, whose distributions form a continuous
      strong Markov process on $H$.
    \item For any two solutions $X$ and $Y$ of (\ref{1.2}) we have for
      all $t\geq s\geq 0$
      \beq\label{1.5}
        \beg{split}
          \|X_t - Y_t\|_H^2
          &\le \big\{\|X_s-Y_s\|_H^{1-r} + (r-1)(\eta  -\theta)
            \ll_1^{(r+1)/2} (t-s) \big\}^{-2/(r-1)}\\
          &\le \|X_s-Y_s\|_H^2\land \big\{(r-1)(\eta
            -\theta)\ll_1^{(r+1)/2} (t-s) \big\}^{-2/(r-1)}.
        \end{split}
      \end{equation}
      Consequently, setting $P_tF(x):= \mathbb EF(X_t(x))$
      for $F:H\to\R$, Borel measurable, so that the expectation makes sense, we
      have
      that $(P_t)_{t>0}$ is a Feller semigroup on $C_b(H)$ and, in addition,
      for Lipschitz continuous $F$
      \beq\label{eq1.9.1}
        \big| P_t F(x) - P_t F(y) \big|
        \leq \scr L (F) \| x-y \|_H,\ \ \ x,y\in H,
      \end{equation} where $\scr L(F)$ is the Lipschitz constant
      of $F$.
    \item $P_t$  has a unique invariant probability
      measure $\mu$ and for some constant $C>0$, $\mu$ satisfies
      \beq\label{1.6}
        \sup_{x\in H} |P_t F(x)-\mu(F)|\le C \scr L (F) t^{-1/(r-1)},
        \quad t>0,
      \end{equation}
      for any Lipschitz continuous function $F$ on
      $H$.
      Moreover, $\mu(\|\cdot\|_{r+1}^{r+1})<\infty.$
    \item If $\si> \dd$ then for any two solutions $X$ and $Y$ of
      (\ref{1.2}) we have for all $t\geq s\geq 0$
      \beq\label{eq1.10.1}
        \| X_t - Y_t \|_H \leq \|X_s-Y_s\|_H \, e^{-(\sigma-\delta)(t-s)}
      \end{equation}
      and there exists $C>0$ such that
      \beq\label{1.10}
        \| X_t-Y_t \|_H \le C \e^{-(\si -\dd)t},
        \quad t\ge 1.
      \end{equation}
      Consequently, for some constant $C>0$,
      \beq\label{1.11}
        \sup_{x\in H} |P_t F(x)-\mu(F)| \le C\scr L(F) \e^{-(\si -\dd)t},
        \quad t\ge 1,
      \end{equation}
      for any Lipschitz continuous function $F$ on $H$.
  \end{enumerate}
\end{thm}

\begin{rem}\label{rem1.2}
  \begin{enumerate}[(1)]
    \item When $Q=0$,  the Dirac measure $\dd_0$ is the unique invariant
      measure. Thus, (\ref{1.6}) with $F(x):=\|x\|_H$ implies
      \[
        \sup_{x}\|X_t(x)\|_H \le C t^{-1/(r-1)},
        \quad t>0.
      \]
      This coincides with the optimal decay of the solution to the
      classical porous medium equation obtained by Aronson and Peletier (see
      \cite[Theorem 2]{AP}).
    \item In the case where $\Phi=0$ and $\Psi(r)= \aa r +r^m$ for
      $\aa\ge 0$ and $m\ge 3$ odd, and $L:=\DD$ on a regular
      domain in $\R^d$, in \cite{BBDR} and \cite{DR1} much stronger integrability
      results for the invariant measure have been proved,
      namely, if either $m=3$ or $\aa>0$ then
      $\mu\big(\big| \nabla(\mathop{\text{sign}} x \, |x|^n) \big|^2\big)<\infty$
      for any $n\ge 1.$
    \item In the case where $L:=\DD$ on a bounded smooth domain in $\R^d$,
      the existence of an invariant measure $\mu$
      was proved in \cite{BDR} under
      the conditions that $\kk_0 |s|^{r-1} \le \Psi'(s)\le C \kk_1 |s|^{r-1}$ and
      $|\Phi(s)|\le C +\dd |s|^r$ for some constants
      $C, \kk_0,\kk_1>0, \dd\in (0, 4\kk_0 \ll_1 (r+1)^{-2}$ and all $s\in\R.$
      Also in \cite{BDR} stronger integrability properties for $\mu$ have
      been proved, namely that
      $\mu\big(\big| \nabla(\mathop{\text{sign}} x \, |x|^\ell) \big|^2\big)<\infty$
      for all $\ell\in[(r+1)/2, r]$.
  \end{enumerate}
\end{rem}
Finally, we note that in this paper the coefficient in front of the noise
is constant (i.e.\ so-called additive noise). Under the usual Lipschitz
assumptions, however, properly reformulated versions of our results also
hold for non-constant diffusion coefficients. 
Details on this will be contained in a forthcoming paper.

\section{Some preliminaries}

We shall make use of a finite-dimensional approximation argument
to construct the solution of (\ref{1.2}).
For any $n\ge 1$, let $r_t^{(n)}:=(r_{t,1}^{(n)},\cdots, r_{t,n}^{(n)})$
solve the following SDE on $\R^n$:
\beq\label{2.1}
  \d r_{t,i}^{(n)} = q_i \d B_t^i -\ll_i
  {\bf m}\Big(e_i \Psi\Big(\sum_{k=1}^n r_{t, k}^{(n)}e_k\Big)\Big)
  \d t+{\bf m}\Big(e_i \Phi\Big(\sum_{k=1}^n r_{t,
  k}^{(n)}e_k\Big)\Big)\d t
\end{equation}
with $r_{0, i}^{(n)}= \<X_0,e_i\>, \ 1\le i\le n$,
where $X_0 : \Omega\to H$ is a fixed $\F_0 / \scr B (H)$-measurable
map such that $\EE \|X_0\|_H^2 < \infty$.
Here and below for a topological space $S$ we denote its Borel
$\sigma$-algebra by $\scr B(S)$.
By \cite[Theorem~1.2]{K} there exists a unique solution to
(\ref{2.1}) for all $t\geq 0$.

\beg{lem}\label{L2.1}
  Under the assumptions of Theorem \ref{T1.1}, there
  exists a constant $C>0$ independent of $n$ and $X_0$ such that
  $X_t^{(n)}:= \sum_{i=1}^n r_{t,i}^{(n)}e_i$ satisfies
  \beq\label{2.2}
    \mathbb E\int_0^T {\bf m}(|X_t^{(n)}|^{r+1})\d t\le C(\|X_0\|_H^2 +T),
    \ \ T>0,
  \end{equation}
  and
  \beq\label{2.3}
    \mathbb E \|X_t^{(n)}\|_H^2 \le C(1+ t^{-2/(r-1)}),
    \ \ \ t>0.
  \end{equation}
\end{lem}

\beg{proof}
  By  (\ref{1.3}) we have
  \beq\label{NN}
    s\Psi(s)
    \ge s\Psi(0) +\eta |s|^{r+1} +\si s^2,
    \ \ s,t\in\R,
  \end{equation}
  and by (\ref{1.3'}) we have
  \[
    -\m \big( \Phi(x) \, L^{-1}x \big)
    \leq -\Phi(0) \, \m(L^{-1}x) + \theta\|x\|_{r+1}^{r+1} +\delta\|x\|_2^2 ,
    \quad x\in L^{r+1}(\m).
  \]
  Hence for all $x\in \mathop{\text{span}}\{e_i : i\in\N\}$
  \begin{multline*}
    -\m \bigl(\Psi(x)\, x\bigr)-\m\bigl( \Phi(x) \, L^{-1}x \bigr)\\
    \leq \bigl(|\Psi(0)| + |\Phi(0)|\,\lambda_1^{-1}\bigr)\|x\|_2
      -(\eta - \theta)\|x\|_{r+1}^{r+1} - (\sigma-\delta)\|x\|_2^2  .
  \end{multline*}
  Combining this with
  (\ref{2.1}) and using It\^o's formula, we obtain
  \beg{equation}\label{NN2}
    \frac 12\d  \|X_t^{(n)}\|_H^2
    \le\d M_t^{(n)} + \frac{c_1}{2} \d t - \frac{c_2}{2} {\bf m}(|X_t^{(n)}|^{r+1})\d t
  \end{equation}
  for some local martingale $M_t^{(n)}$ and constants $c_1,c_2>0$ independent of $n$.
  This implies (\ref{2.2}).
  Moreover, since
  ${\bf m}(|X_t^{(n)}|^{r+1}) \ge \ll_1^{(r+1)/2}\|X_t^{(n)}\|_H^{r+1},$
  it follows  from (\ref{NN2}) that
  \beq\label{**}
    \EE \|X_t^{(n)}\|_H^2 - \EE\|X_s^{(n)}\|_H^2
    \leq c_1(t-s) - c_2 \int_s^t \bigl( \EE\|X_u^{(n)}\|_H^2 \bigr)^{\frac{r+1}{2}} \,\d u,
    \quad 0\leq s\leq t.
  \end{equation}
  To prove (\ref{2.3}), let $h$ solve the equation
  \beq\label{A1}
    h'(t)=-c_2 h(t)^{(r+1)/2} + c_1,
    \quad t\geq 0, \qquad
    h(0)= \EE \|X_0\|_H^2 + (4 c_1/c_2)^{2/(r+1)}.
  \end{equation}
  Then it is
  easy to see that  (\ref{**}) implies
  \beq\label{A2}
    \mathbb E\|X_t^{(n)}\|_H^2\le h(t),\ \ \ t\ge
    0.
  \end{equation}
  Let
  $\phi_t:= h(t)- \mathbb E\|X_t^{(n)}\|_H^2$ and
  \[
    \tau:=\inf\{t\ge 0: \phi_t\le 0\}.
  \]
  Suppose $\tau < \infty$, then by continuity $\phi_\tau \leq 0$
  and by the mean-value theorem and (\ref{**}), (\ref{A1}) we obtain
  \[
    \phi_t
    \ge \phi_0 - c_2 \int_0^t \Big(h_\vv(u)^{(r+1) / 2}- \bigl(\mathbb E\|X_u^{(n)}\|_H^2\bigr)^{(r+1) / 2}\Big)\d u
    \ge (4 c_1 / c_2)^{2 / (r+1)} -c\int_0^t \phi_u\d u,
    \ \ 0\le t\le \tau,
  \]
  where $c:= c_2  \frac{r+1}{2} \max_{t\in [0,\tau]} t^{(r-1)/2} = c_2 \frac{r+1}{2} \tau^{(r-1)/2}$.
  By Gronwall's lemma we arrive at
  $\phi_\tau \ge (4 c_1 / c_2)^{2/(r+1)} \e^{-c\tau}>0$.
  This contradiction proves (\ref{A2}).

  To estimate $h(t),$ let
  \[
    \tau:= \inf\{t\ge 0: h(t)^{(r+1)/2}\le 2c_1/c_2\}.
  \]
  Since $h(0)^{(r+1)/2}\ge 4c_1/c_2> 2c_1/c_2,\ \tau\ge t_0$
  for some $t_0>0$ independent of $n$. Indeed, we may define
  $t_0$ as $\tau$ above with $h$ replaced by the solution to (\ref{A1})
  with initial condition
  $h(0) :=(4c_1/c_2)^{2/(r+1)}.$ By (\ref{A1}) we have
  \[
    h'(t) \le -\ff {c_2} 2 h(t)^{(r+1)/2},
    \ \ \ \ 0 \le t\le \tau.
  \]
  Therefore, for some constant $c>0$ independent of $n$,
  \beq\label{A3}
    h(t)\le ct^{-2/(r-1)},
    \ \ \ 0\le t\le \tau .
  \end{equation}
  Clearly,
  $h'(t)\le 0$ for
  all $t\geq 0$, since by an elementary consideration we have
  $h \ge (c_1/c_2)^{2/(r+1)}$, consequently
  \[
    h(t)\le h(\tau)\le c\tau^{-2/(r-1)}\le ct_0^{-2/(r-1)},
    \ \ \ t>\tau.
  \]
  Therefore, (\ref{2.3}) holds.
\end{proof}

According to (\ref{2.2}) in Lemma~\ref{L2.1}, $X^{(n)}$ is bounded
in $L^{r+1}(\R_+\times\OO\times E; \nu\times P\times {\bf m})$,
where $\nu(\d t):= \e^{-t}\d t.$ Thus, there exists a subsequence
$n_k\to \infty$ and a process $X$ such that $X^{(n_k)}\to X$
weakly in $L^{r+1}(\R_+\times\OO\times E; \nu\times P\times {\bf
m})$. To prove that this limit provides a solution of (\ref{1.2}),
we shall make use of Theorem 3.2 in Chapter 1 of \cite{KR}. We
state this result in detail for the reader's convenience
specialized to our situation.


\beg{thm}\label{T2.2}
  {\bf (\cite[Theorem I.3.2]{KR})}
  Consider three maps $v: \R_+\times\Omega\to\BB$,
  $\tilde v : \R_+ \times \Omega\to\BB'$, $h:\R_+\times \Omega\to H$
  such that
  \begin{enumerate}[(i)]
    \item $v$ is $\scr B(\R_+)\otimes \F / \scr B(\BB)$-measurable and
      $v_t := v(t,\,\cdot\,)$ is $\F_t / \scr B(\BB)$-measurable
      for all $t\geq 0$.
    \item $\tilde v$ is $\scr B(\R_+)\otimes \F / \scr B(\BB')$-measurable
      and $\tilde v_t := \tilde v (t,\,\cdot\,)$ is
      $\F_t /\scr B(\BB')$-measurable. Moreover,
      $\int_0^T \|\tilde v_t\|_{\BB'}\,\d t <\infty$ $P$-a.s.\ for all $T>0$.
    \item $h$ is an $H$-valued $(\F_t)$-adapted continuous local semi-martingale.
  \end{enumerate}
  Set
  \[
    \tt h_t:= \int_0^t \tt v_s\d s +h_t.
  \]
  If $\tt h_t(\oo)= v_t(\oo)$ for $\nu\times P$-a.e.\ $(t,\oo)$,
  then $\tt h_t$ is an
  $H$-valued continuous $(\F_t)$-adapted process satisfying the following
  It\^o formula for the square of the norm:
  \beq\label{2.8}
    \|\tt h_t\|_H^2
      = \|\tt h(0)\|_H^2 + 2 \int_0^t {}_{\BB'}\<\tilde v_s , v_s \>_{\BB} \d s
        + 2 \int_0^t \<\tt h_s, \d h_s\>_H +[h]_t.
  \end{equation}
  where $[h]$ denotes the quadratic variation process of $h$.
\end{thm}

\section{Proof of the existence}
\paragraph{a)}
By Lemma \ref{L2.1} and (\ref{1.3}), $\{\Psi(X_t^{(n)})\}$ and
$\{\Phi(X_t^{(n)})\}$ are   bounded in $L^{(r+1)/r}(\R_+\times
\OO\times E; \nu\times P\times {\bf m})$, where $\nu(\d t):=
\e^{-t}\d t.$ Hence there exist a subsequence $n_k\to \infty$  and
processes $U,V\in L^{(r+1)/r}(\R_+\times \OO\times E; \nu\times
P\times {\bf m})$ such that \beq\label{3.1}
  \Psi(X^{(n_k)}) \to  U,\ \Phi(X^{(n_k)})\to V\
  \text{weakly\ in}\ L^{(r+1)/r}(\R_+\times \OO\times E; \nu\times P\times {\bf m}).
\end{equation}
Moreover, by Lemma \ref{L2.1}, we may also assume that
\beq\label{3.1'}
  X^{(n_k)}\to \bar X\
  \text{weakly\ in}\  L^{r+1}(\R_+\times \OO\times E; \nu\times P\times {\bf m}).
\end{equation}
E.g.\ by \cite[Chap.~3, \S{}~7]{diestel} we may also assume that
the Cesaro means of the sequences in (\ref{3.1}) converge strongly
in $L^{(r+1)/r}(\R_+ \times \Omega \times E ; \; \nu\times P\times
\m )$ so the limits have $\scr B(\R_+) \otimes \F \otimes \scr
M$-measurable versions. Furthermore, as continuous processes the
approximants are all progressively measurable as $L^{(r+1) /
r}(\m)$-valued processes, hence so are  their limits. In
particular, these are adapted. The same holds for the sequence in
(\ref{3.1'}) respectively its limit with $(r+1)/r$ replaced by
$r+1$. Below we always consider versions of $U$, $V$, $\bar X$
with all these measurability properties and denote them by the
same symbols. Since for $t\geq 0$
\[
  {\bf m}(X_t^{(n_k)}e_i)
  = \<X_0,e_i\>
    + \int_0^t \big\{ {\bf m}(\Psi(X_s^{(n_k)}) L e_i)
                      + {\bf m}(e_i \Phi(X_s^{(n_k)})
               \big\}
      \d s
    + q_i B_t^i,
    \ \ \ 1\le i\le n_k,
\]
it follows from (\ref{3.1}) and (\ref{3.1'})  that, for any real-valued
bounded measurable process $\varphi$,
\[
  \mathbb E\int_0^T \varphi_t {\bf m}(\bar X_t e_i) \nu(\d t)
  = \mathbb E\int_0^T \varphi_t
    \bigg\{ \<X_0, e_i\>
            + \int_0^t \big\{{\bf m}(U_s L e_i)+{\bf m}(V_s e_i)\big\}\d s
            + q_i B_t^i
    \bigg\}
    \nu(\d t)
\]
for all $T>0$. Thus,
\beq\label{3.2}
  {\bf m}(\bar X_t e_i)  =
  \<X_0, e_i\> +\int_0^t  {\bf m}(U_s L e_i + V_s e_i)\d s+ q_i B_t^i,
  \quad \text{for $\nu\times P$-a.e.\ } (t,\oo), i\ge 1.
\end{equation}

\paragraph{b)}
To apply Theorem \ref{T2.2}, let
\[
  \tilde v_s:= \bar L U_s + J V_s .
\]
By (\ref{1.3}), Lemma~\ref{L2.1}, Proposition~\ref{prp1.1}
and Corollary~\ref{cor1.1}
we have $\mathbb E \int_0^T \|\bar v_s\|_{\B'}\d s<\infty$ for any $T>0$.
So, we see
that in Theorem~\ref{T2.2} conditions
(i), (ii) with $v=\bar X$
and also (iii) with $h := QW$
are satisfied and by Proposition~\ref{prp1.1} and Corollary~\ref{cor1.1},
(\ref{3.2}) with $e_i$ replaced by $\lambda_i^{-1}e_i$ implies
\[
  {}_{\BB'}\<\bar X_t, e_i\>_{\BB}
  = {}_{\BB'}\<X_0, e_i\>_{\BB}
    + \int_0^t {}_{\BB'}\<\tilde v_s, e_i\>_{\BB} \;\d s
    + {}_{\BB'}\< QW_t , e_i \>_{\BB} ,
  \quad i\geq 1, \; \text{for $\nu\times P$-a.e.\ $(t,\omega)$.}
\]
Hence defining
\beq\label{3.4}
  X_t
  := X_0 + \int_0^t \tt v_s \d s + QW_t,
  \quad t\geq 0
\end{equation}
we see that
\begin{equation}\label{eq3.5.1}
  \bar X = X \quad \text{$\nu\times P$-a.e.}
\end{equation}
Therefore,
by Theorem~\ref{T2.2}, $X$ is an $H$-valued continuous
$(\scr F_t)$-adapted process
and (\ref{2.8}) holds with $X$ replacing $\tilde h$.
Therefore, to prove that $X$ solves (\ref{1.2}),
by Proposition~\ref{prp1.1} and Corollary~\ref{cor1.1} it suffices to show that
\beq\label{3.5}
  {\bf m}\big(e_i[V_s -\ll_i U_s]\big)
  = {\bf m}\big(e_i [\Phi(\bar X_s)-\ll_i \Psi(\bar X_s)]\big),
  \quad i\ge 1, \; \text{for $\nu\times P$-a.e.\ $(s,\omega)$.}
\end{equation}
This will  be proved by the following two steps.

\paragraph{c)}
We claim that for any $\psi : \R_+ \to \R_+$ bounded, Borel measurable
with compact support,
\beq\label{3.6}
  \liminf_{k\to\infty} \int_0^\infty \psi(t)\mathbb E \|X_t^{(n_k)}\|_H^2 \d t
  \ge \int_0^\infty \psi(t)\mathbb E\|\bar X_t\|_H^2\d t .
\end{equation}
Since $X^{(n_k)}\to \bar X$ weakly in $L^2(\R_+\times \OO\times E;
\nu\times  P\times  {\bf m})$, by Fatou's lemma we have
\beg{equation*}
  \beg{split}
    \int_0^\infty\psi(t)\mathbb E \|\bar X_t\|_H^2\d t
    &= \sum_{i=1}^\infty \ll_i^{-1}\mathbb E
       \int_0^\infty\psi(t)  {\bf m}(e_i\bar X_t)^2 \d t\\
    &= \sum_{i=1}^\infty \ll_i^{-1}\lim_{k\to\infty}\mathbb E
       \int_0^\infty \psi(t) {\bf m}(e_iX_t^{(n_k)}) {\bf m}(e_i\bar X_t)\d t\\
    &\le \ff 1 2 \sum_{i=1}^\infty \ll_i^{-1}\liminf_{k\to\infty}\mathbb E
       \int_0^\infty \psi(t)  {\bf m}(e_iX_t^{(n_k)})^2 \d t\\
    &\qquad +  \ff 1 2 \sum_{i=1}^\infty \ll_i^{-1}\mathbb E
       \int_0^\infty \psi(t)  {\bf m}(e_i\bar X_t)^2 \d t\\
    &\le \ff 1 2 \liminf_{k\to\infty}\mathbb E
       \int_0^\infty \psi(t) \|X_t^{(n_k)}\|_H^2 \d t +  \ff 1 2 \mathbb E
       \int_0^\infty \psi(t)  \|\bar X_t\|_H^2 \d t.
  \end{split}
\end{equation*}
Since $\bar X\in L^2(\R_+\times \OO\times E; \nu\times P\times {\bf
m})$ so that $\int_0^\infty\psi(t) \mathbb E \|\bar X_t\|_H^2\d
t<\infty,$ this implies (\ref{3.6}) immediately.

\paragraph{d)}
By \eqref{2.8} in Theorem~\ref{T2.2}, \eqref{eq1.5} and \eqref{eq1.7.1}
we have
\beq\label{3.8}
  \mathbb E\|X_t\|_H^2
  = \EE \|X_0\|_H^2 -2\int_0^t \mathbb E \big({\bf m}(\bar X_s U_s)
    + \m (L^{-1}(\bar X_s) V_s) \big) \d s+\sum_{i=1}^\infty \ll_i^{-1} q_i^2 t.
\end{equation}
On the other hand, by It\^o's formula,
\beq\label{2.4}
  \mathbb E\|X_t^{(n)}\|_H^2
  = \EE \|X_0\|_H^2 - 2\mathbb E\int_0^t
        \m \bigl( X_s^{(n)}\Psi(X_s^{(n)})
                  + L^{-1}(X_s^{(n)}) \Phi(X_s^{(n)})
           \bigr)
      \d s
    +\sum_{i=1}^n \ll_i^{-1} q_i^2t.
\end{equation}
Then for any $\varphi\in L^{r+1}(\R_+\times \OO\times E; \nu\times
P\times{\bf m})$, we obtain from (\ref{1.3}), (\ref{1.3'}) and
(\ref{2.4}) that \beg{equation}
  \beg{split}\label{FF}
    0\le I_k(t)
    &:= 2\mathbb E \int_0^t
        \m \Bigl( [X_s^{(n_k)}-\varphi_s] \bigl[\Psi(X_s^{(n_k)})- \Psi(\varphi_s)\bigr] \\
    &\qquad \qquad
        + L^{-1} (X_s^{(n_k)}-\varphi_s) \bigl[ \Phi(X_s^{(n_k)}) - \Phi(\varphi_s)\bigr]\Bigr)
        \d s\\
    &\hphantom{:}= 2 \mathbb E \int_0^t
          \m \Bigl( X_s^{(n_k)} \Psi(X_s^{(n_k)}) + L^{-1} (X_s^{(n_k)}) \Phi(X_s^{(n_k)})\Bigr)
        \d s\\
    &\quad
        - 2 \mathbb E\int_0^t \m \Bigl(
          \varphi_s \bigl[ \Psi(X_s^{(n_k)})-\Psi(\varphi_s) \bigr]
          + X_s^{(n_k)} \Psi(\varphi_s) \\
    &\qquad\qquad
          + L^{-1} (X_s^{(n_k)})\Phi(\varphi_s)
          + L^{-1} (\varphi_s) \bigl[ \Phi(X_s^{(n_k)}) - \Phi(\varphi_s) \bigr] \Bigr)
        \d s\\
    &\hphantom{:}= - \mathbb E \|X_t^{(n_k)}\|_H^2
       + \|X_0\|_H^2
       + \sum_{i=1}^{n_k} \ll_i^{-1} q_i^2 t \\
    &\quad
       - 2 \mathbb E\int_0^t
       \m
       \Bigl(
           \varphi_s \bigl[ \Psi(X_s^{(n_k)})-\Psi(\varphi_s) \bigr]
           + X_s^{(n_k)}\Psi(\varphi_s) \\
    &\qquad\qquad
           + L^{-1} (X_s^{(n_k)}) \Phi(\varphi_s)
           + L^{-1} (\varphi_s) \bigl[ \Phi(X_s^{(n_k)}) - \Phi(\varphi_s) \bigr]
       \Bigr)
       \d s,
  \end{split}
\end{equation}
Since $\Psi(X^{(n_k)})\to U$ and $\Phi(X^{(n_k)})\to V$ weakly in
$L^{(r+1)/r}(\R_+\times \OO\times E; \nu\times P\times{\bf m})$ and
$X^{(n_k)}\to\bar X$ weakly in $L^{r+1}(\R_+\times \Omega\times E,
\; \nu\times P\times \m)$, for any $\psi : \R_+\to\R_+$, bounded,
Borel-measurable with compact support, we obtain from (\ref{FF})
and (\ref{3.6}) that \beg{equation*}
  \beg{split}
    0 \le
    &\int_0^\infty \psi(t)\d t
     \Big\{
       -\mathbb E\|\bar X_t\|_H^2 +\EE\|X_0\|_H^2+\sum_{i=1}^\infty
       \ll_i^{-1} q_i^2 t\\
    &  - 2 \mathbb E\int_0^t
       \m \Bigl(
         \varphi_s \bigl[ U_s-\Psi(\varphi_s) \bigr]
         + \bar X_s\Psi(\varphi_s)
         + L^{-1} (\bar X_s) \Phi(\varphi_s)
         + L^{-1} (\varphi_s) \bigl[ V_s- \Phi(\varphi_s) \bigr]
       \Bigr)
       \d s
     \Big\}.
  \end{split}
\end{equation*}
Combining this with (\ref{3.8}) we arrive at
\beq\label{II}
  0\le \int_0^\infty\psi(t)\d t
  \mathbb E \int_0^t \m
  \Bigl(
    [ \bar X_s-\varphi_s ] \bigl[ U_s-\Psi(\varphi_s) \bigr]
    + L^{-1} (\bar X_s-\varphi_s) \bigl[ V_s-\Phi(\varphi_s) \bigr]
  \Bigr)
  \d s.
\end{equation}
By first taking $\varphi_s:=\bar X_s-\vv\tt\varphi_s e_i$ for
given $\vv>0$ and $\tt\varphi\in L^{\infty}(\R_+\times \OO;
\nu\times P)$, then dividing by $\vv$ and letting $\vv\to 0$, we
obtain from the continuity of $\Phi,\Psi$ and the dominated
convergence theorem, which is valid due to (\ref{1.3}),
(\ref{1.3'}) and because $\bar X\in L^{r+1}([0,\infty)\times
\OO\times E;\nu\times P\times {\bf m})$, that
\[
  \int_0^\infty \psi(t)\d t
  \mathbb E\int_0^t\tt\varphi_s{\bf m}
  \big(
    e_i[U_s-\Psi(\bar X_s)-\ll_i^{-1}( V_s-\Phi(\bar X_s))]
  \big)
  \d s\ge 0.
\]
Since $\psi : \R_+\to\R_+$ bounded, Borel-measurable with compact
support, and $\tt\varphi\in L^{\infty}(\R_+\times \OO; \nu\times
P)$ are arbitrary, this implies (\ref{3.5}).

\section{Proof of the other assertions}
\paragraph{a)}
Proofs of the uniqueness, \eqref{eq1.8.1}, (\ref{1.5}) and (\ref{1.6}).\\
\eqref{eq1.8.1} is an immediate consequence of \eqref{2.3}
and \eqref{3.6}, since $X$ is $P$-a.e.\ continuous in $t$.
Let $X$ and $Y$ be solutions of (\ref{1.2}). Then by \eqref{eq1.8}
$Z:= X-Y$ solves the equation
\[
  Z_t = Z_0 + \int_0^t \Bigl[ \bar L \bigl( \Psi(X_t)-\Psi(Y_t) \bigr)
                               + J \bigl( \Phi(X_t)-\Phi(Y_t) \bigr) \Bigr]\d t,
  \quad t\geq 0.
\]
Thus, by Theorem~\ref{T2.2} with $h=0$,
\eqref{eq1.5}, \eqref{eq1.7.1} and finally by
(\ref{1.3}), (\ref{1.3'}),
we obtain
\beg{equation}\label{U}
  \|Z_t\|_H^2
  = \|Z_0\|_H^2 - 2\int_0^t
      \m \Bigl(
        [X_s-Y_s] \bigl[ \Psi(X_s)-\Psi(Y_s) \bigr]
        + L^{-1} (X_s-Y_s) \bigl[ \Phi(X_s)-\Phi(Y_s) \bigr]
      \Bigr)
    \d s
  \le 0.
\end{equation}
If $X_0=Y_0$,
this implies  $Z_t=0$ for all $t\geq 0$.
By a slight modification of a
standard argument one obtains as usual that the uniqueness also
implies the stated Markov property and hence the
semigroup property of $P_t$.
The strong Markov property then follows from the Feller
property of $P_t$ proved below, since all solutions of \eqref{1.2}
have continuous sample paths in $H$.

Similarly to (\ref{U}), we have for $0\leq s\leq t$ that
\beq\label{3.9}
  \beg{split}
    &\|X_t-Y_t\|_H^2 - \|X_s-Y_s\|_H^2 \\
    &\le - 2\int_s^t
      \big\{
        (\si- \dd)\|X_u-Y_u\|_2^2 +(\eta
        -\theta) \|X_u-Y_u\|_{r+1}^{r+1}
      \big\}
      \d u.
  \end{split}
\end{equation}
Noting that $\si\ge \dd,\eta>\theta$,
$\|\,\cdot\,\|_2^2 \geq \lambda_1 \|\,\cdot\,\|_H^2$,
$\|\cdot\|_{r+1}^{r+1}\ge \ll_1^{(r+1)/2}\|\cdot\|_H^{r+1}$ and
that for $\vv>0$ the function
$h_{\vv,t} :=  \{(\vv+\|X_s-Y_s\|_H)^{1-r} + (r-1)
(\eta-\theta)\ll_1^{(r+1)/2} (t-s)\}^{-2/(r-1)}$, $t\geq s$, solves
for $s\geq 0$ fixed the equation
\beq\label{eq4.2.1}
  h_t'= -2 (\eta -\theta)\ll_1^{(r+1)/2} h_t^{(r+1)/2},
  \quad t\geq s(\geq 0),
  \quad h_0= (\|X_s-Y_s\|_H +\vv)^2
\end{equation}
due to the same comparison argument as in the proof of Lemma~\ref{L2.1},
it follows that $\|X_t-Y_t\|_H^2 \leq h_{\vv,t}$ $\forall\; t\geq s$.
Letting $\vv\to 0$ this
implies (\ref{1.5})
and the
Feller property of $P_t$.
By \eqref{eq1.8.1} it follows that $P_t |F|(x) <\infty$ for all Lipschitz
continuous  $F:H\to\R$. Now \eqref{eq1.9.1} is obvious by \eqref{1.5}.

\paragraph{b)}
Proof of (3).\\
Let $\dd_0$ be the Dirac measure at $0\in H$. Set
\[
  \mu_n
  := \ff 1 n \int_0^n \dd_0 P_t \d t,
  \ \ \ n\ge 1.
\]
We intend to show the tightness of $\{\mu_n\}$.
Then by the Feller property of $P_t$ the weak
limit of a subsequence provides an invariant probability measure
of $P_t$. By Lemma \ref{L2.1} and the weak convergence of
$X^{(n_k)}$ to $X$, we have
\[
  \int_H \m(|x|^2)\mu_n(\d x) = \ff 1 n \int_0^n \mathbb E
  \m(|X_t(0)|^2)\d t\le C
\]
for some constant $C>0$ and all $n\ge 1,$ where we set
$\m(|x|^2)=\infty$ if $x\notin L^2(\m)$. Since
the function $x\mapsto \m(|x|^2)$ is compact, that is,
$\{x\in H: \m(|x|^2)\le r\}$ is relatively compact in
$H$ for any $r\ge 0$, we conclude that $\{\mu_n\}$ is tight.

Next, let $\mu$ be an invariant probability measure.
For any bounded Lipschiz function $F$ on $H$, (\ref{1.5}) implies
that there exists $C>0$ such that
\[
  |P_t F(x)-\mu(F)|
  \le \int_H\mathbb E|F(X_t(x))-F(X_t(y))|\mu(\d y)\le C\scr L(F) t^{-1/(r-1)}
\]
for all $x\in H, t>0.$
Thus, (\ref{1.6}) holds and hence $P_t$ has a unique invariant measure.

Let $\mu$ be the invariant probability measure of $P_t$.
It remains to show that $\mu(\|\cdot\|_{r+1}^{r+1})<\infty$.\\
By \eqref{eq1.8.1}, since $\mu$ is $P_t$-invariant, we have
\beg{equation}\label{LLL}
    \int_H\|x\|_H^2\mu(\d x)
    = \int \EE \bigl\|X_1(x)\bigr\|_H^2 \; \mu(\d x) \leq 2C < \infty.
\end{equation}
Next, since $X$ is the weak limit of $X^{(n_k)}$ in
$L^{r+1}([0,\infty)\times \OO\times E; \nu\times P\times {\bf m})$,
by H\"older's inequality we have \beg{equation*}
  \beg{split}
    & \int_1^2 \mathbb E\,\m(|X_t(x)|^{r+1})\d t =\lim_{k\to\infty} \int_1^2
      \mathbb E\m \big(X_t^{(n_k)} (|X_t|^r \text{sgn}(X_t))\big) \d t\\
    & \le \liminf_{k\to\infty}
      \bigg(
        \int_1^2 \mathbb E \m (|X_t^{(n_k)}|^{r+1})\d t
      \bigg)^{1/(r+1)}
      \bigg(
        \int_1^2 \mathbb E \m(|X_t|^{r+1})\d t
      \bigg)^{r/(r+1)}.
  \end{split}
\end{equation*}
Therefore, by \eqref{2.2}, for some $C_1>0$
\[
  \int_1^2 \mathbb E\,\m(|X_t(x)|^{r+1})\d t \le C_1(1+\|x\|_H^2),
  \quad x\in H.
\]
Then by (\ref{LLL}) we obtain
\[
  \mu(\|\cdot\|_{r+1}^{r+1})
  = \int_H\mu(\d x)\int_1^2 \mathbb E\,\m(|X_t(x)|^{r+1})\d t
  < \infty.
\]

\paragraph{c)}
Exponential ergodicity, i.e.\ proof of (4).\\
If $\si>\dd$ then  (\ref{3.9}) implies
\[
  \|X_t - Y_t\|_H^2
  \le \|X_s-Y_s\|_H^2 \, \e^{-2(\si -\dd)(t-s)},
  \quad t\ge s\ge 0.
\]
So, \eqref{eq1.10.1} holds.
Combining this with (\ref{1.5}) we arrive at
\[
  \|X_{t+1}-Y_{t+1}\|_H^2
  \le \|X_1 - Y_1\|_H^2 \, \e^{-2(\si -\dd)t}\le C\e^{-2(\si-\dd)t},
  \quad t\ge 0.
\]
This implies (\ref{1.10}) and
(\ref{1.11}) immediately.
\qed

\beg{thebibliography}{99}

\bibitem{Aronson}
  D.G. Aronson,
  \emph{The porous medium equation,}
  Lecture Notes Math. Vol. 1224, Springer, Berlin, 1--46, 1986.

\bibitem{AP}
  D.G. Aronson and L.A. Peletier,
  \emph{Large time behaviour of solutions of the
     porous medium equation in bounded domains,}
  J. Diff. Equ. 39(1981), 378--412.

\bibitem{BBDR}
  V. Barbu, V.I. Bogachev,  G. Da Prato and M. R\"ockner,
  \emph{Weak solution to the stochastic porous medium  equations: the degenerate case,}
  preprint.

\bibitem{BDR}
  V.I. Bogachev,  G. Da Prato and M. R\"ockner,
  \emph{Invariant measures of   stochastic generalized porous medium equations,}
  to appear in Dokl. Math.

\bibitem{DR1}
  G. Da Prato and M. R\"ockner,
  \emph{Weak solutions to stochastic porous media equations,}
  J. Evolution Equ. 4(2004), 249--271.

\bibitem{DR2}
  G. Da Prato and M. R\"ockner,
  \emph{Invariant measures for a stochastic porous medium equation,}
  preprint SNS, 2003; to appear in Proceedings of
  Conference in Honour of K. It\^o, Kyoto, 2002.

\bibitem{diestel}
  J. Diestel,
  \emph{Geometry of Banach spaces -- selected topics},
  Lect. Notes Math. 485,
  Springer, 1975.

\bibitem{K}
  N.V. Krylov,
  \emph{On Kolmogorov's equations for finite-dimensional diffusions,}
  Lecture Notes Math. 1715, 1--63, Springer, Berlin 1999.

\bibitem{KR1}
  N.V. Krylov and M. R\"ockner,
  \emph{Strong solutions of stochastic equations with singular time dependent drift,}
  to appear in Probab. Theory Relat. Fields.

\bibitem{KR}
  N.V. Krylov and B.L. Rozovskii,
  \emph{Stochastic evolution equations,}
  Translated from Itogi Naukii Tekhniki, Seriya Sovremennye Problemy Matematiki
  14(1979), 71--146, Plenum Publishing Corp. 1981.

\bibitem{ma-roe92}
  Z.M. Ma and M. R\"ockner,
  \emph{Introduction to the theory of (non-symmetric) Dirichlet forms},
  Springer, 1992.
\end{thebibliography}

\end{document}